\documentclass[10pt]{amsart}
\usepackage[mathscr]{eucal}
\usepackage{amssymb, amsmath,array, amscd}
\usepackage{enumerate}

\theoremstyle{change}
\newtheorem{thm}{Theorem}[section]
\newtheorem{THM}{Theorem}
\newtheorem{prop}{Proposition}[section]

\newtheorem{remark}{Remark}[section]
\newtheorem{cor}{Corollary}[section]
\newtheorem{conjecture}{Conjecture}[section]
\newtheorem{problem}{Problem}[section]
\newtheorem{example}{Example}[section]

\begin{document}

\title[Planar Webs]{On Planar Webs with Infinitesimal Automorphisms}
\author{D. Mar{\'\i}n}
\address{Departament de Matem{\`a}tiques \\ Universitat Aut{\`o}noma de Barcelona \\ E-08193  Bellaterra (Barcelona)\\ Spain}
\email{davidmp@mat.uab.es}

\author{J. V. Pereira}
\address{Instituto de Matem{\'a}tica Pura e Aplicada\\ Est.
D. Castorina, 110\\
22460-320, Rio de Janeiro, RJ, Brasil}
\email{jvp@impa.br}

\author{L. Pirio}
\address{IRMAR \\ Campus de Beaulieu \\ 35042 Rennes Cedex, France}
\email{luc.pirio@univ-rennes1.fr  }

\thanks{The second author is supported by Cnpq and Instituto Unibanco. The third
author was partially supported by the International Cooperation Agreement Brazil-France.}

\begin{abstract} We investigate the space of abelian relations 
of planar webs admitting  infinitesimal automorphisms. As an application 
 we construct $4k - 14$ new  algebraic families of global exceptional $k$-webs on the 
projective plane, for each $k\ge 5$. 
\end{abstract}

\maketitle

\date{\today}

\section{Introduction and statement of the results} 

\subsection{Planar Webs}
A germ of regular $k$-web $\mathcal W=\mathcal F_1 \boxtimes \cdots \boxtimes \mathcal F_k$ 
on $(\mathbb C^2,0)$ is a collection of 
$k$ germs of smooth foliations $\mathcal F_i$ subjected to the condition that 
any two of these foliations have distinct tangent spaces at the origin.

One of the most intriguing invariants of a web is its {\it space of 
abelian relations} $\mathcal A(\mathcal W)$.
If the foliations $\mathcal F_i$ are induced by $1$-forms $\omega_i$ then by definition
\begin{align*}
\mathcal A(\mathcal W) = \Big\lbrace {\big(\eta_i\big)}_{i=1}^k \in
(\Omega^1(\mathbb C^2,0))^k \, \, \Big| \, \,\forall i \, \,
d\eta_i =0 \,, \; \eta_i\wedge \omega_i=0 \, \, \text{ and } \sum_{i=1}^k\eta_i = 0 \Big\rbrace
\, .
\end{align*}
The dimension of $\mathcal A(\mathcal W)$ is commonly called the
{\it rank} of $\mathcal W$ and noted by 
$\mathrm{rk}(\mathcal W)$. It is a theorem of Bol that $\mathcal A(\mathcal W)$ is
a finite-dimensional $\mathbb C$-vector space and moreover
\begin{align}
\label{bolbound}
\mathrm{rk}( \mathcal W) \le \frac{1}{2}\, (k-1)(k-2)\, .
\end{align}

An interesting chapter of the theory of webs concerns the characterization
of webs of {\it maximal rank}, {\it i.e} webs for which (\ref{bolbound}) is in fact an equality. It follows from Abel's Addition Theorem that all
the webs $\mathcal W_C$ obtained from reduced plane curves $C$
  by projective duality  are of maximal rank ({\it cf.} \S \ref{S:action} for 
details). The webs analytically equivalent to
some $\mathcal W_C$ are the so called {\it algebrizable webs}.

It can be traced back to Lie a remarkable result  that 
says that all $4$-webs of maximal rank are in fact algebrizable.
In the early 1930's Blaschke claimed to have extended Lie's result
to $5$-webs of maximal rank. Not much latter Bol came up with 
a counter-example: a $5$-web of maximal rank that is not algebrizable.

The non-algebrizable 
webs of maximal rank are nowadays called  {\it exceptional webs}. For a long time Bol's web 
remained as the only example of exceptional planar web in the literature. 
The following
quote  illustrates quite well this fact.

\begin{quote}
{\it (\ldots)  we cannot refrain from mentioning what we consider to be
the fundamental problem on the subject, which is to determine the maximum
rank non-linearizable webs. The strong conditions must imply that there are
not many. It may not be unreasonable to compare the situation with the
exceptional simple Lie groups.}
\begin{flushright}
Chern and Griffiths in \cite{Jbr}.
\end{flushright}
\end{quote}

A comprehensive account of the  current state of the art
concerning the exceptional webs  is available at \cite[Introduction \S 3.2.1]{theseluc}, \cite{Robert} and \cite[\S 1.4]{PT}. Here we will just mention
that before this work no exceptional $k$-web with $k\ge 10$ appeared
in the literature. 

At first glance, the list of known exceptional webs up today does not
reveal common features among them. Although at a second look one sees that 
many of them (but not all,  not even the majority) have one property in common:
infinitesimal automorphisms.

\subsection{Infinitesimal Automorphisms}
In \cite{cartan}, {\'E}. Cartan proves that 
{\it a 3-web which admits an 2-dimensional continuous group of transformations is hexagonal}.
It is then an exercise to deduce  that a $k$-web
 ($k> 3$) which admits  $2$ linearly independent infinitesimal automorphisms is parallelizable and
in particular  algebrizable.

Cartan's result naturally leads to the following question:

\begin{quote}
{\it What can be said about  webs 
which admit one infinitesimal  automorphism?}  
\end{quote}

In fact, Cartan  answers this question for 3-webs. In {\it loc. cit.} 
he establishes that such a web is equivalent to those induced by the $1$-forms  $dx,dy, dy -u(x+y)dx$, where
$u$ is a germ of holomorphic function.

It is very surprising that this story stops here\ldots \,
To our knowledge, there is no other study concerning webs with
infinitesimal automorphisms,  
although they are particularly interesting. Indeed, on the one hand
their study is considerably simplified by 
the presence of an infinitesimal automorphism, but on the other hand, 
these webs can be very interesting from a geometrical point of view:  we will
show they are connected to the theory of exceptional webs.

\subsection{Variation of the Rank} 
Let ${\mathcal W}$ be a regular web in $({\mathbb C}^2,0)$
which admits an infinitesimal automorphism $X$, {\it i.e.} $X$ 
is a germ of vector field whose local flow preserves the foliations
of $\mathcal W$. As we will see in \S\ref{S:geral} 
the Lie derivative $L_X=i_Xd + di_X$ with respect to $X$  
induces a linear operator  on $\mathcal A(\mathcal W)$. 
Most of our results  will follow 
from an analysis of such operator.

In \S\ref{S:liouville}  we use this 
operator to give  a simple description of the  abelian relations of $\mathcal W$ 
and from this  we will deduce  in \S\ref{S:rank} what we 
consider~our~main~result:

\begin{THM}\label{T:1}
Let ${\mathcal W}$ be a $k$--web which admits a transverse infinitesimal automorphism $X$. Then
\[
\mathrm{rk}(\mathcal W \boxtimes \mathcal F_X) =\mathrm{rk}(\mathcal W) + (k -1)\, .
\]
In particular,  ${\mathcal W}$ is of maximal rank if and only if  
$\mathcal W \boxtimes \mathcal F_{X}$ is of maximal rank.
\end{THM}

We will derive from Theorem \ref{T:1} the existence of 
new families of exceptional webs.
\subsection{New Families of Exceptional Webs}

If we start with a reduced plane curve $C$  invariant
under an algebraic ${\mathbb C}^*$-action on ${\mathbb P}^2$ then we
obtain  a dual algebraic ${\mathbb C}^*$-action on 
$\check{\mathbb P}^2$, 
letting invariant the algebraic web ${\mathcal W}_C$ ({\it cf.} \S\ref{S:action} for
details). Combining this construction  with Theorem \ref{T:1} 
we deduce~our~second~main~result. 

\begin{THM}\label{T:2}
For every $k  \geq 5$   
there exist a family of dimension at least  $\lfloor k/2 \rfloor -1$ of 
pairwise non-equivalent  exceptional global $k$-webs on ${\mathbb P}^2$.
\end{THM}

In fact, for each $k\ge 5$,  we obtain  $4k - 15$ other  families
of smaller dimension.

We also give a complete classification of all the 
exceptional 5-webs of the type $ {\mathcal W}\boxtimes {\mathcal F}_X$ 
where $X$ is an infinitesimal automorphism of $ {\mathcal W}$ ({\it cf.} Corollary \ref{C:classi}).


\section{Generalities on webs with infinitesimal automorphisms}\label{S:geral}

Let $\mathcal F$ be a regular foliation on $(\mathbb C^2,0)$
induced by a (germ of) $1$-form $\omega$. We
say that a (germ of)  vector field $X$ is an infinitesimal
automorphism of $\mathcal F$ if the foliation $\mathcal F$
is preserved by the local flow of $X$. In algebraic terms: 
$ L_X \omega \wedge \omega = 0 \, .$

When the infinitesimal automorphism $X$ is transverse to $\mathcal F$, 
{\it i.e} when $\omega(X)\neq0$,  
then a simple computation ({\it cf}. \cite[Corollary 2]{Percy}) shows that the
$1$-form
\[
   \eta = \frac{\omega}{i_X \omega}
\]
is   closed and satisfies $L_X \eta =0$. By definition, the integral 
\[
u(z) = \int^z_0 \eta 
\]
is the {\it canonical first integral} of ${\mathcal F}$ (with respect to $X$). Clearly, we have $u(0)=0$ and $L_X(u)= 1$.\vspace{0.15cm}

Now let $\mathcal W$ be a germ of regular $k$-web on $(\mathbb C^2,0)$
induced by the (germs of) $1$-forms $\omega_1, \ldots, \omega_k$ and let 
 $X$ be an  infinitesimal
automorphism of $\mathcal W$. Here, of course, we mean that $X$ is an 
infinitesimal automorphism for~all~the~foliations~in~$\mathcal W$.
 
By hypothesis, we have $L_X\, \omega_i \wedge \omega_i=0$ for $i=1,\ldots,k$. Then because the Lie derivative $L_X$  is linear and commutes with 
$d$, it induces a linear map
\begin{eqnarray}
\label{themap!}
L_X : \mathcal A(\mathcal W) & \longrightarrow & \mathcal A(\mathcal W) \\
(\eta_1,\ldots ,\eta_k) & \longmapsto& (L_X\eta_1 ,\ldots, L_X \eta_k)
\, .  \nonumber
\end{eqnarray}

This map is central in this paper: all our results come from an analysis of 
the $L_X$-invariant subspaces of $ {\mathcal A}(\mathcal W)$.

\section{Abelian relations of webs with infinitesimal automorphisms}

\subsection{Description of ${\mathcal A}({\mathcal W})$ in presence of an infinitesimal automorphism}\label{S:liouville} 
In this section, ${\mathcal W}=\mathcal F_1 \boxtimes \cdots \boxtimes \mathcal F_k $ 
denotes a $k$-web in $({\mathbb C}^2,0)$ which admits 
an infinitesimal automorphism $X$, regular and transverse to the foliations ${\mathcal F}_i$
in a neighborhood of the origin. 

Let $i\in \{1,\ldots,k\}$  be fixed.   
We note ${\mathcal A}^i({\mathcal W})$ 
the vector subspace of $\Omega^1(\mathbb C^2,0)$ spanned 
by the $i$-th components $\alpha_i$ of 
abelian relations $(\alpha_1,\ldots,\alpha_k)\in {\mathcal A}({\mathcal W})$. 
If $u_i=\int \eta_i$ denotes the canonical first integral of ${\mathcal F}_i$ 
with respect to $X$, then for $\alpha_i \in {\mathcal A}^i({\mathcal W})$, 
there exists a holomorphic germ $f_i \in \mathbb C\{t\}$ 
such that $\alpha_i=f_i(u_i)\,du_i$. 

Assume now that ${\mathcal A}^i({\mathcal W})$ is not trivial and 
let $\big\{\alpha_i^\nu=f_\nu(u_i)\,du_i \, | \, \nu=1,\ldots,n_i\}$ be a basis.
Since $L_X : \mathcal A^i(\mathcal W) \to \mathcal A^i(\mathcal W)$ 
is a linear map, there exist complex constants $c_{\nu \mu}$ such that, for $\nu=1,\ldots,n_i$ we have 
\begin{eqnarray}
\label{baz}
L_X(\alpha_i^\nu)=\sum_{\mu=1}^{n_i} c_{\nu \mu}\, \alpha_i^\mu \; .
\end{eqnarray}

But $
L_X(\alpha_i^\nu)=L_X\big(f_\nu(u_i)\,du_i\big)=X\big( 
f_\nu(u_i)\big)du_i+f_\nu(u_i)\,L_X\big(du_i\big)=f_\nu'(u_i)\,{du_i} 
$ for any $\nu$, so relations (\ref{baz}) are equivalent to the scalar ones
\begin{eqnarray}
\label{diffeq}
\quad \quad
f_\nu'
=\sum_{\mu=1}^{n_i} c_{\nu \mu}\, {f_\mu} \, , \quad \quad \, \nu=1,\ldots,n_i \, .
\end{eqnarray}

Now let $\lambda_1,\ldots,\lambda_\tau \in {\mathbb C}$ be the 
eigenvalues of the map $L_X$ acting on ${\mathcal A}({\mathcal W})$ 
corresponding to maximal eigenspaces with corresponding dimensions 
 $\sigma_1,\ldots,\sigma_\tau $. The differential equations 
 (\ref{diffeq}) give us the following description of $ {\mathcal A}({\mathcal W})$:

\begin{prop} \label{description}
The abelian relations of $ {\mathcal W}$ are of the form 
 $$ P_1(u_1)\,e^{\lambda_i\,u_1}\,du_1+\cdots+ P_k(u_k)\,e^{\lambda_i\,u_k}\,du_{k}=0 
$$
where $P_1, \ldots, P_k$ are  polynomials of degree less   or equal to $ \sigma_i$. 
\end{prop}

A non-zero eigenvalue $\lambda$ of the map (\ref{themap!}) corresponds to a functional equation of the form 
$ c_1\,e^{\lambda\,u_1}+\cdots+ c_k\,e^{\lambda\,u_k}={cst.}$ 
where the $c_i$'s are complex constants.
Using Abel's method to solve functional equations (see \cite[Chapter 2]{theseluc}), 
such $\lambda$ can be computed effectively from the $u_i$'s, 
which can be effectively computed from the $\omega_i$'s. 
In a few words:  Proposition \ref{description} gives an effective tool
to compute the abelian relations of $\mathcal W$.


\subsection{Proof of Theorem \ref{T:1}}\label{S:rank}
With Proposition \ref{description} at hand we are able to prove our main result.

Let  $\mathcal W= \mathcal F_1 \boxtimes \cdots \boxtimes \mathcal F_k$
and for $i= 1 \ldots k$, set  $\eta_i = du_i$ as the differential 
of the canonical first integral of $\mathcal F_i$ relatively to $X$. We note $x$ a first integral of the foliation ${\mathcal F}_X$, normalized such that $x(0)=0$.

When $j$ varies from $2$ to $k$, we have
\[
i_X(\eta_1 - \eta_j) =0 \qquad\quad \text{ and } \qquad\quad  L_X(\eta_1 - \eta_j) =
 0 \, . 
\]
Consequently there exists $g_j \in \mathbb C\{x\}$ such that
\begin{equation}\label{E:xxx}
du_1 - du_j - g_j(x)\,dx = 0 \, . 
\end{equation}
Clearly these are abelian relations for the web $\mathcal W \boxtimes \mathcal F_X$.
They span a $(k-1)$-dimensional 
vector subspace $\mathcal V$ of the maximal eigenspace  of $L_X$ associated to the eigenvalue zero, noted 
$\mathcal A_0(\mathcal W\boxtimes \mathcal F_X)$.

Observe that $\mathcal V$ fits in the following exact sequence ($i$ is the natural inclusion): 
\begin{align}
\label{exactseq}
0 \to {\mathcal V} \stackrel{i}{\longrightarrow}
{\mathcal A_0(\mathcal W \boxtimes \mathcal F_X)}\stackrel{L_X}{\longrightarrow}
 {\mathcal A_0(\mathcal W)}   \,.
\end{align}
Indeed, the kernel  $ K:= \ker \{ L_X : \mathcal A_0(\mathcal W \boxtimes \mathcal F_X) 
\to
 {\mathcal A_0(\mathcal W)} \}$ is generated by abelian relations of the form
$\sum_{i=1}^k c_i du_i + g(x)\,dx = 0$, where $c_i \in \mathbb C$ 
and $g \in \mathbb C\{x\}$.   Since $i_X du_i=1$ for each $i$,
 it follows that the constants
$c_i$   satisfy $\sum_{i=1}^k c_i =0$. It implies that 
 the abelian relations in the kernel
of $L_X$ can be written as linear combinations of abelian relations 
of the form (\ref{E:xxx}).  Therefore 
\begin{equation}\label{E:yyy}
K =  \mathcal V
\end{equation}
and consequently $\ker L_X \subset \mathrm{Im } \, \, i$. The exactness
of (\ref{exactseq}) follows easily.   

From general principles we deduce that the sequence
\[
0 \to \frac{\mathcal V}{\mathcal A_0(\mathcal W)\cap \mathcal V} \stackrel{i}{\longrightarrow}
\frac{\mathcal A_0(\mathcal W \boxtimes \mathcal F_X)}{\mathcal A_0(\mathcal W)} \stackrel{L_X}{\longrightarrow}
 \frac{\mathcal A_0(\mathcal W)}{L_X \mathcal A_0(\mathcal W)} \, ,
\]
is also exact. Thus to  prove the Theorem it suffices to 
verify the following assertions:
\begin{enumerate}
\item[(a)] $\mathcal V$ is isomorphic to $$
\frac{\mathcal V}{\mathcal A_0(\mathcal W)\cap \mathcal V} \oplus 
  \frac{\mathcal A_0(\mathcal W)}{L_X \mathcal A_0(\mathcal W)} \, ;
$$
\item[(b)] the morphism $L_X:\mathcal A_0(\mathcal W \boxtimes \mathcal F_X) \to
 \mathcal A_0(\mathcal W)$
 is surjective;
 \item[(c)] the vector spaces 
 $$\frac{\mathcal A_0(\mathcal W \boxtimes \mathcal F_X)}{\mathcal A_0(\mathcal W)} 
 \, \, \text{ and }  \, \,
 \frac{\mathcal A(\mathcal W \boxtimes \mathcal F_X)}{\mathcal A(\mathcal W)}$$ are isomorphic.
\end{enumerate}

To verify assertion (a), notice that the nilpotence of $L_X$ on $\mathcal A_0(\mathcal W)$
implies that $\frac{\mathcal A_0(\mathcal W)}{L_X \mathcal A_0(\mathcal W)}$ is isomorphic 
to $\mathcal A_0(\mathcal W) \cap K$. Combined 
with (\ref{E:yyy}), it implies assertion (a).\vspace{0.1cm}

To prove assertion (b), it suffices to construct a map $\Phi:\mathcal A_0(\mathcal W)
\to \mathcal A_0(\mathcal W \boxtimes \mathcal F_X)$ such that $ L_X \circ \Phi = \mathrm{Id}$.
Proposition \ref{description}  implies that $\mathcal A_0(\mathcal W)$ is spanned by abelian 
relations of the form
$ \sum_{i=1}^k c_i u_i^r du_i =0  ,$
where $c_i$ are complex numbers and $r$ is a non-negative integer. 
For such an abelian relation, since 
\[
   \sum_{i=1}^k c_i u_i^r du_i = \frac{1}{r+1}\,L_X \Big( \sum_{i=1}^k c_i u_i^{r+1} du_i \Big) =0 \, ,
\]
there exists an unique $g\in \mathbb C \{ x\}$ satisfying
$  \sum_{i=1}^k c_i u_i^{r+1} du_i + g(x)\, dx =0 \, .
$
If we set 
\[
\Phi\Big(\sum_{i=1}^k c_i u_i^r du_i \Big) 
= \frac{1}{r+1}\, \Big(  \sum_{i=1}^k c_i u_i^{r+1} du_i + g(x)\, dx \Big) \, 
\]
then  $L_X \circ \Phi= \mathrm{Id}$ on $\mathcal A_0(\mathcal W)$ and assertion (b) follows.

To prove assertion (c) we first notice that 
\[
\mathcal A(\mathcal W \boxtimes \mathcal F_X)=\mathcal A_0(\mathcal W \boxtimes \mathcal F_X)  
\oplus
  \mathcal A_*(\mathcal W \boxtimes \mathcal F_X) \, 
\]
where $\mathcal A_*(\mathcal W \boxtimes \mathcal F_X)$ denotes the sum of eigenspaces corresponding to non-zero eigenvalues.
Of course $\mathcal A_*(\mathcal W \boxtimes \mathcal F_X)$ is invariant and moreover we have the equality
\[
L_X\big(\mathcal A_*(\mathcal W \boxtimes \mathcal F_X)\big) = \mathcal A_*(\mathcal W \boxtimes \mathcal F_X) \, .
\]
But $L_X$ {\it kills} the $\mathcal F_X$-components of abelian relations. 
In particular, it implies  
\[
L_X\big(\mathcal A_*(\mathcal W \boxtimes \mathcal F_X)\big) \subset  \mathcal A_*(\mathcal W ).
\]
This is sufficient to show that $\mathcal A_*(\mathcal W \boxtimes \mathcal F_X) = \mathcal A_*(\mathcal W )$ and
deduce assertion (c) and, consequently that
\[
\mathrm{rk}(\mathcal W \boxtimes \mathcal F_X) =\mathrm{rk}(\mathcal W) + (k -1)\, .
\]

Because $k(k-1)/2=(k-1)(k-2)/2+(k-1)$, the above inequality  implies immediately the last assertion of Theorem 1. 
\qed


\section{New Families of exceptional webs}

\subsection{Algebrizable Webs with Infinitesimal Automorphisms}\label{S:action}

Let $C \subset \mathbb P^2$ be a degree $k$ reduced curve.
 If $U\subset \check{\mathbb P}^2$ is a simply-connected open set not intersecting
 $\check C$ and if 
$\gamma_1, \ldots, \gamma_k: U  \to C$
are the holomorphic maps defined by the intersections of lines in $U$ with $C$
then Abel's Theorem implies that 
\[
\mathrm{Tr}(\omega)= \sum_{i=1}^k \gamma_i^* \omega = 0 \,
\]
for every $\omega \in H^0(C,\omega_C)$, where $\omega_C$ denotes the dualizing sheaf of $C$.

 The maps $\gamma_i$
define  the $k$-web $\mathcal W_C$ on $U$ and the trace formula above associates
an abelian relation of $\mathcal W_C$ to each $\omega \in H^0(C,\omega_C)$. 
 Since $h^0(C,\omega_C)=(k-1)(k-2)/2$, the web $\mathcal W_C$ is of maximal rank.

Suppose now that  $C$ is invariant by a $\mathbb C^*$-action
$\varphi: \mathbb C^* \times \mathbb P^2 \to \mathbb P^2.$
Notice that  $\varphi$ induces a dual action
 $\check{\varphi}:\mathbb C^* \times
\check{\mathbb P^2} \to \check{\mathbb P^2}$ satisfying
$
\varphi_t \circ \gamma_i = \gamma_i \circ \check{\varphi_t}
$ for $i=1,\ldots,k$. 
Consequently the web $\mathcal W_C$ admits an infinitesimal automorphism.

In a suitable projective coordinate system $[x:y:z]$, a plane curve $C$
invariant by a ${\mathbb C}^*$-action is cut out by an equation of the form
\begin{equation}\label{E:aluffi}
x^{\epsilon_1}\cdot y^{\epsilon_2}\cdot z^{\epsilon_3} \cdot \prod_{i=1}^k (  x^a + \lambda_i y^b z^{a-b}  ) \, 
\end{equation}
where $\epsilon_1, \epsilon_2, \epsilon_3 \in \{ 0,1 \}$,  $k,a,b \in \mathbb N$ are such that $k\ge1$, $a\ge 2$, $1\le b\le a/2$,
$\gcd(a,b)=1$ 
and the $\lambda_i$ are distinct non zero complex numbers ({\it cf.} \cite[\S 1]{aluffi}
for instance). Notice that here the $\mathbb C^*$-action in question is 
\begin{equation}\label{E:zzz}
\begin{array}{l c l c l}
\varphi &:&  \mathbb C^* \times \mathbb P^2  &\to& \mathbb P^2 \\
 & &  ( t, [x:y:z] ) &\mapsto& [t^{b(a-b)}x:t^{a(a-b)} y: t^{ab}z ] \, .
\end{array}
\end{equation}

Moreover once we fix $\epsilon_1,\epsilon_2,\epsilon_3,k,a,b$  we can always choose $\lambda_1=1$ and in this 
case the set of $k-1$ complex numbers $\{\lambda_2,\ldots, \lambda_k\}$ projectively characterizes the 
curve $C$. In particular one promptly sees that 
 there exists a $(d-1)$-dimensional family  of degree 
 $2d$ (or $2d +1$)  
 reduced plane curves all projectively
distinct and invariant by the same  $\mathbb C^*$-action: for a given $2d+ \delta$ 
with $\delta\in \{0,1\}$ set $a=2$, $b=1$, $\epsilon_1=\delta$ and $\epsilon_2=\epsilon_3=0$.

A moment of reflection shows that the number of discrete parameters giving
distinct families of degree $d$ curves of the form (\ref{E:aluffi}) is
\[
\underbrace{\Big\lfloor \frac{d}{2} \Big\rfloor}_{\epsilon_1=\epsilon_2=\epsilon_3=0} 
+ \underbrace{ 3 \Big\lfloor \frac{d-1}{2} \Big\rfloor}_{\epsilon_i=\epsilon_j=0, \,  \epsilon_k=1}
+ \underbrace{3\Big\lfloor \frac{d-2}{2} \Big\rfloor}_{\epsilon_i=\epsilon_j=1, \,  \epsilon_k=0}
 + \underbrace{\Big\lfloor \frac{d-3}{2} \Big\rfloor}_{\epsilon_1=\epsilon_2=\epsilon_3=1} 
 -\; 2  = 4d - 10 \, .
\]

Notice that the $-2$ appears on left hand side because the  curves $\{y=0\}$ and $\{z=0\}$ are indistinguishable
when $a=2$. 

\subsection{Proof of Theorem \ref{T:2}} 
If $C$ is a reduced curve of the form (\ref{E:aluffi}) then $\mathcal W_C$ is invariant
by an algebraic $\mathbb C^*$-action $\check \varphi$. We will note by $X$ the infinitesimal
generator of $\check \varphi$ and by $\mathcal F_X$ the corresponding foliation. 
From the discussion on  the last
paragraph, Theorem \ref{T:2} follows at once from the stronger:

\begin{thm}\label{T:geral}
If $\deg C\ge 4$ then $\mathcal W_C \boxtimes \mathcal F_X$ is exceptional. Moreover if $C'$ 
is another curve invariant by $\varphi$ then $\mathcal W_C \boxtimes \mathcal F_X$  is analytically
equivalent to $\mathcal W_{C'} \boxtimes \mathcal F_X$ if and only if the curve $C$ is projectively equivalent to
$C'$.
\end{thm}
\begin{proof} Since $\mathcal W_C$ has maximal rank it follows from Theorem \ref{T:1}
that $\mathcal W_C \boxtimes \mathcal F_X$ is also of maximal rank. Suppose that its
localization  at a point $p\in \mathbb P^2$ is algebrizable and
let $\psi:(\mathbb P^2,p) \to (\mathbb C^2,0)$ be a holomorphic algebrization. Since 
both $\mathcal W_C$ and $\psi_*(\mathcal W_C)$ are linear webs of maximal rank 
it follows from a result of Nakai \cite{nakai} that $\psi$ is the localization of 
an automorphism of $\mathbb P^2$.  But the generic leaf of $\mathcal F_X$
is not contained in any line of $\mathbb P^2$ and consequently $\psi_*(\mathcal W \boxtimes
\mathcal F_X)$ is not linear. This concludes the proof of the theorem.
\end{proof}

\begin{remark}\rm
We do not know if the families above are {\it irreducible} in the sense
that the generic element does not admit a deformation as an exceptional
web that is not contained in the family. Due to the presence of automorphism
one could imagine that they are indeed degenerations of some other 
exceptional webs.
\end{remark}

\subsection{A characterization  result}
Combining Theorem \ref{T:1} with Lie's Theorem we can easily prove the

\begin{cor}\label{C:classi}
Let $\mathcal W$  be a 4-web that admits a transverse  infinitesimal automorphism $Y$.
If $\mathcal W\boxtimes \mathcal F_Y$ is exceptional then it is analytically equivalent
to an exceptional 5-web $\mathcal W_C \boxtimes \mathcal F_X$ described in Theorem \ref{T:geral}.
\end{cor}
\begin{proof}
It follows from Theorem \ref{T:1} that $\mathcal W$ is of maximal rank. 
Lie's Theorem implies 
that $\mathcal W$ is analytically equivalent to $\mathcal W_C$ for some reduced plane quartic $C$.
Since the local flow of $Y$ preserves $\mathcal W$ there exists a (germ) of vector field $X$
whose local flow preserves $\mathcal W_C$. Using again Nakai's result we deduce that the
germs of automorphisms on the local flow of $X$ are indeed projective automorphisms. This is 
sufficient to conclude that $X$ is a global vector field preserving $\mathcal W_C$.
\end{proof}

The example below shows that
Theorem \ref{T:geral} does not give all the exceptional webs admitting an infinitesimal
automorphism.

\begin{example}\rm
In \cite{crasluc}, it is proved that the web ${\mathcal W}$ 
 induced by the functions $x,y,x+y,x-y,x^2+y^2$ is 
exceptional. Moreover it admits the radial vector field 
$R=x\,\partial/\partial_x+y\,\partial/\partial_y$  
as a transverse infinitesimal automorphism. 
Theorem \ref{T:1} implies that the $6$-web 
$\mathcal W \boxtimes \mathcal F_R$  is also exceptional.
This result was previously obtained by 
determining an explicit basis of the space of abelian 
relations, see \cite[p. 253]{theseluc}.  
\end{example}


\section{Problems}

\subsection{A conjecture about the nature of the abelian relations}
It is clear from Proposition \ref{description} that for 
webs $\mathcal W$ admitting infinitesimal automorphisms there
exists a  Liouvillian extension of the field of definition of $\mathcal W$
containing all its abelian relations. We believe that a similar statement
should hold for arbitrary webs $\mathcal W$.

\begin{conjecture}
The abelian relations  of a web are Liouvillian.
\end{conjecture}

Our belief is supported by the recent works of H{\'e}naut \cite{Henaut} and Ripoll 
\cite{Ripoll} on abelian relations and of Casale \cite{Guy} 
on non-linear differential Galois Theory.

When $\mathcal W$ is of maximal rank the main result of \cite{Henaut} shows that
there exists a Picard-Vessiot extension of the field
of definition of $\mathcal W$ containing all the abelian relations. 
In the general case, one should be able to deduce a similar result from
the above mentioned work of Ripoll.

On the other hand, and at least over polydiscs, \cite[Theorem 6.4]{Guy}  implies 
that the foliations with first integrals on Picard-Vessiot extension are 
transversely projective. Since the first integrals in question are components
of abelian relations they are of finite determinacy and hopefully this should
imply that they are indeed Liouvillian.

\subsection{Restricted Chern's Problem}
With the techniques now available, the classification of all exceptional 5-webs (``{\it Chern's problem}'') seems completely out of reach. 
So we propose the

\begin{problem}
Classify exceptional 5-webs admitting infinitesimal automorphisms.
\end{problem} 

Notice that this restricted version is not completely hopeless. The linear map $L_X$  can be ``integrated''
giving birth to a holomorphic  action on $\mathbb P(\mathcal A(\mathcal W))$.  The Poincar{\'e}-Blaschke
curves will be orbits of this action and the dual action
will induce an automorphism of the associated Blaschke 
surface. This seems valuable extra data that may lead to a solution of
the restricted Chern's  problem. 

For a definition of the above mentioned concepts 
 see \cite[Chapter 8]{theseluc}.


\end{document}